\documentstyle[11pt,A4]{article}
\begin{document}

\def \a{\alpha}
\def \b{\beta}
\def \g{\gamma}
\def \d{\delta}
\def \l{\lambda}
\def \m{\mu}
\def \n{\nu}
\def \r{\rho}
\def \s{\sigma}
\def \tr{\mathrm{tr}}

\begin{titlepage}

\title{ Heisenberg Double and Pentagon Relation}

\author{R.M. Kashaev\thanks{On leave from
St. Petersburg Branch of the Steklov Mathematical Institute,
Fontanka 27, St. Petersburg 191011, RUSSIA}\\ \\
Laboratoire de Physique Th\'eorique
{\sc enslapp}\thanks{URA 14-36 du CNRS,
associ\'ee \`a l'E.N.S. de Lyon,
et \`a l'Universit\`e de Savoie}
\\
ENSLyon,
46 All\'ee d'Italie,\\
69007 Lyon, FRANCE\\
E-mail: {\sf rkashaev\@@enslapp.ens-lyon.fr}}

\date{March 1995}

\maketitle

\abstract{It is shown that the Heisenberg double has a canonical element,
satisfying the pentagon relation. From a given invertible constant solution
to the pentagon relation one can restore the structure
of the underlying algebras. Drinfeld double can be realized as
 a subalgebra in the tensor square of the Heisenberg double. This enables
 one to write down solutions to the Yang-Baxter relation in terms of solutions
to the pentagon relation.
}
\vskip 2cm

\rightline{{\small E}N{\large S}{\Large L}{\large A}P{\small P}-L-512/95}

\end{titlepage}

\section{Introduction}

The theory of quantum groups \cite{D} appeared as an algebraic setting for the
 construction of solutions to the Yang-Baxter equation \cite{Y,B}:
\begin{equation}
R_{12}R_{13}R_{23}=R_{23}R_{13}R_{12}.   \label{eq1.1}
\end{equation}
According to \cite{D}, for any Hopf algebra ${\cal A}$ one can construct a
 quasi-triangular Hopf algebra $D({\cal A})$, called Drinfeld double,
 where
the universal $R$-matrix, satisfying relation (\ref{eq1.1}), exists.
 In \cite{FRT}
 it was shown how to reconstruct the underlying Hopf algebra from a given
 solution to equation (\ref{eq1.1}). In \cite{RS,AF,S} another, Heisenberg
 double $H({\cal A})$, has been introduced, which, unlike the Drinfeld
double, is not a Hopf algebra.

 The purpose of the present letter is to show that the constant ``pentagon''
 relation
\begin{equation}
S_{12}S_{13}S_{23}=S_{23}S_{12},   \label{eq1.2}
\end{equation}
plays the same role in Heisenberg double as Yang-Baxter relation (\ref{eq1.1})
 does
in Drinfeld's one. This is the content of Sect.~\ref{sec2}. In Sect.~\ref{sec3}
we show that the Drinfeld double can be realized as
 a subalgebra in the tensor square of the Heisenberg double. This gives an
explicit formula for solutions to the Yang-Baxter equation in terms of
 solutions to the pentagon equation. As an example, in Sect.~\ref{sec4} the
Borel subalgebra
 of $U_q(sl(2))$ at $|q|<1$ is considered and a slightly generalized form
of the quantum dilogarithm identity of \cite{FK} is obtained.

{\bf Acknowledgements.}
It is a pleasure to thank L.D. Faddeev, J.M. Maillet and
M.A. Semenov-Tian-Shansky for the valuable discussions.

\section{Heisenberg Double and Pentagon Relation}

\label{sec2}
\setcounter{equation}{0}
Let $\{e_\a\}$ be a linear basis of associative and co-associative bialgebra
${\cal A}$, with the following multiplication and co-multiplication rules:
\begin{equation}
e_\a e_\b=m_{\a \b}^\g e_\g,\quad \Delta(e_\a)=\mu_\a^{\b\g}e_\b\otimes
 e_\g,\label{eq2.1}
\end{equation}
where summations over the repeated indices are implied.
The Heisenberg double
$H({\cal A})$ can be defined as an associative algebra with generating
elements $\{e^\b,e_\a\}$ and the following defining relations:
\begin{equation}
e_\a e_\b=m_{\a \b}^\g e_\g,\quad e^\a e^\b=\mu^{\a \b}_\g e^\g,\quad
e_\a e^\b=m^\b_{\r\g}\mu^{\g\s}_\a e^\r e_\s.     \label{eq2.2}
\end{equation}
There are two subalgebras of $H({\cal A})$ with the linear bases
$\{e_\a\}$ and $\{e^\a\}$, which are equivalent to the algebra ${\cal A}$
 and it's dual ${\cal A^*}$. The analog of the adjoint representation is given
by a realization in terms of the structure constants:
\begin{equation}
\langle\a|e_\b|\g\rangle=m_{\a\b}^\g,\quad
\langle\a|e^\b|\g\rangle=\mu^{\b\g}_\a.        \label{eqadj}
\end{equation}
The remarkable property of the Heisenberg double is described by the following
 theorem.

\noindent
\newtheorem{theorem}{Theorem}
\begin{theorem}The canonical element $S=e_\a\otimes e^\a$ in the
 Heisenberg double satisfies the pentagon relation (\ref{eq1.2}).
\end{theorem}

Now, following \cite{FRT}, for a given invertible solution $S$ of the pentagon
 relation (\ref{eq1.2}), define two bialgebras ${\cal B}$ and ${\cal B^*}$,
 generated with entries of matrices $F$ and $G$, respectively, and the
 following multiplication and co-multiplication rules:
\begin{equation}
F_1F_2S_{12}=S_{12}F_1,\quad S_{12}G_1G_2=G_2S_{12},     \label{eq2.3}
\end{equation}
\begin{equation}
\Delta(F_1)=F_1\otimes F_1,\quad \Delta(G_1)=G_1\otimes G_1.\label{eq2.4}
\end{equation}
Note, that the associativity conditions are ensured by the pentagon relation
(\ref{eq1.2}). In fact, these
algebras are dual as bialgebras with the following pairing:
\begin{equation}
\langle G_1,F_2\rangle=S_{12}.           \label{eq2.5}
\end{equation}
To extract the structure constants for given linear bases
$\{e_\a\}$ in ${\cal B}$ and $\{e^\a\}$ in ${\cal B^*}$ with the canonical
 pairing:
\begin{equation}
\langle e^\b,e_\a\rangle=\delta_\a^\b,        \label{eq2.6}
\end{equation}
represent matrices $F$ and $G$ as linear combinations of the base elements with
some coefficients:
\begin{equation}
F_1=F_1^\a e_\a,\quad G_1=G_{1,\a}e^\a.   \label{eq2.7}
\end{equation}
Substituting these expressions into (\ref{eq2.5}), and using (\ref{eq2.6}),
we obtain:
\begin{equation}
S_{12}=G_{1,\a}F_2^\a.                      \label{eq2.8}
\end{equation}
Expanding the known matrix $S_{12}$ in the form of (\ref{eq2.8}), one can
 calculate
the matrices $G_{1,\a}$ and $F_1^\a$ up to an invertible trasformation over
 the index $\a$. From (\ref{eq2.8}) it follows also, that for a finite
dimensional
 matrix $S$ one associates in this way a finite dimensional algebra with
 the dimension, given by the formula:
\begin{equation}
\mathrm{dim}({\cal B})=\mathrm{rank}(P_{12}S_{12})^{t_1},    \label{eq2.9}
\end{equation}
where $P_{12}$ is the permutation matrix, and the upper index $t_1$ means the
 transposition in the first subspace. Introduce now the dual matrices
 $F_{1,\a}$ and $G_1^\a$ through the equations:
\begin{equation}
\tr_1(F_{1,\a}F_1^\b)=\delta_\a^\b,\quad \tr_1(G_{1,\a}G_1^\b)=\delta_\a^\b.
   \label{eq2.10}
\end{equation}
Using these, one can derive the following formulas for the structure constants:
\begin{equation}
m_{\a\b}^\g=\tr_1(G_{1,\a}G_{1,\b}G_1^\g),\quad
\mu^{\a\b}_\g=\tr_1(F_1^\a F_1^\b F_{1,\g}).    \label{eq2.11}
\end{equation}
 There is an analog of the adjoint representation for these algebras, given by
 the $S$-matrix itself:
\begin{equation}
F_1=S_{01},\quad G_1=S_{10} ,                     \label{eq2.12}
\end{equation}
where the $0$-th subspace corresponds to the representation space. In fact,
 formulas (\ref{eq2.12}) realize the Heisenberg double, defined by
(\ref{eq2.3}) and the
mixed permutation relation of the form:
\begin{equation}
G_1S_{12}F_2=F_2G_1.            \label{eq2.13}
\end{equation}

Formulas (\ref{eq2.9})--(\ref{eq2.11}) make sense for finite dimensional case,
 while the
 infinite dimensional case needs a further qualification. Our result can be
 stated as the following theorem.

\noindent
\begin{theorem}For any invertible  solution to the
 constant pentagon relation (\ref{eq1.2}) one can associate a pair of mutually
 dual associative and co-associative bialgebras.
\end{theorem}

\section{Yang-Baxter and Pentagon Relations}
\label{sec3}
\setcounter{equation}{0}
Let $H({\cal A})$ be the Heisenberg double, defined in Section~\ref{sec2}.
 Co-multiplications of subalgebras ${\cal A}$ and ${\cal A^*}$
can not be extended to any co-multiplication of the whole algebra.
This is the main difference between Heisenberg and Drinfeld doubles.
 It is possible, however, to realize Drinfeld double as a subalgebra
 in the tensor product of two Heisenberg's,
 $H({\cal A})\otimes \tilde H({\cal A})$, where the second ``tilded''
double is defined as follows:
\begin{equation}
\tilde e_\a\tilde e_\b=m_{\a\b}^\g\tilde e_\g,\quad
\tilde e^\a\tilde e^\b=\mu^{\a\b}_\g\tilde e^\g,\quad
\tilde e^\b\tilde e_\a=\mu_\a^{\s\g}m_{\g\r}^\b\tilde e_\s\tilde e^\r.
                                                         \label{eq3.1}
\end{equation}
The canonical element $\tilde S=\tilde e_\a\otimes\tilde e^\a$, satisfies
 the ``reversed'' pentagon relation:
\begin{equation}
\tilde S_{12}\tilde S_{23}=\tilde S_{23}\tilde S_{13}\tilde S_{12}.
\label{eq3.2}
\end{equation}
Using relations (\ref{eq2.2}) and (\ref{eq3.1}) it is easy to show that
 the elements
\begin{equation}
E_\a=\mu_\a^{\b\g} e_\b\otimes\tilde e_\g,\quad
E^\a=m^\a_{\g\b} e^\b\otimes\tilde e^\g.\label{eq3.3}
\end{equation}
 satisfy the defining relations of the Drinfeld double:
\begin{equation}
E_\a E_\b=m_{\a \b}^\g E_\g,\quad E^\a E^\b=\mu^{\a \b}_\g E^\g,\quad
\mu_\a^{\s\g} m^\b_{\g\r}E_\s E^\r=m^\b_{\r\g}\mu^{\g\s}_\a E^\r E_\s.
                                             \label{eq3.4}
\end{equation}
In particular, for the canonical element $R=E_\a\otimes E^\a$ one gets a
 factorized formula:
\begin{equation}
R_{12,34}=S''_{14}S_{13}\tilde S_{24}S'_{23}, \label{eq3.5}
\end{equation}
where
\begin{equation}
S'=\tilde e_\a\otimes e^\a,\quad S''=e_\a\otimes\tilde e^\a. \label{eq3.6}
\end{equation}
By construction, $R$-matrix (\ref{eq3.5}) satisfies the Yang-Baxter relation:
\begin{equation}
R_{1\tilde1,2\tilde2}R_{1\tilde1,3\tilde3}R_{2\tilde2,3\tilde3}=
R_{2\tilde2,3\tilde3}R_{1\tilde1,3\tilde3}R_{1\tilde1,2\tilde2}. \label{eq3.7}
\end{equation}
In fact, it is a consequence of eight different pentagon relations,
 two homogeneous ones, (\ref{eq1.2}) and (\ref{eq3.2}), and six mixed ones:
\begin{eqnarray}
S'_{12}S'_{13}S_{23}=S_{23}S'_{12}, & &
\tilde S_{12}S'_{23}=S'_{23}S'_{13}\tilde S_{12},\nonumber \\
S_{12}S''_{13}S''_{23}=S''_{23}S_{12}, & &
S''_{12}\tilde S_{23}=\tilde S_{23}S''_{13}S''_{12}, \nonumber \\
S'_{12}\tilde S_{13}S''_{23}=S''_{23}S'_{12}, & &
S''_{12}S'_{23}=S'_{23}S_{13}S''_{12}.\label{eq3.8}
\end{eqnarray}
Consider now the case, where algebra ${\cal A}$ is a Hopf
algebra.
The unit element and co-unity map have the form:
\begin{equation}
1=\varepsilon^\a e_\a,\quad \epsilon(e_\a)=\varepsilon_\a,    \label{eq3.9}
\end{equation}
and there are also the antipode and it's inverse maps:
\begin{equation}
\gamma(e_\a)=\gamma_\a^\b e_\b,\quad \overline\gamma(e_\a)=
\overline \gamma_\a^\b e_\b,\quad
\gamma_\a^\g\overline\gamma_\g^\b=\delta_\a^\b .     \label{eq3.10}
\end{equation}
In this case the second double $\tilde H({\cal A})$ can be realized through
the first one $H({\cal A})$:
\begin{equation}
\tilde e_\a=\gamma_\a^\b \overline e_\b,
\quad \tilde e^\a=\overline\gamma^\a_\b \overline e^\b,\label{eq3.11}
\end{equation}
where the overline means the opposite multiplication, which can be realized
by the transposition operation. The four different $S$-matrices
can be expressed now only in terms of the original one:
\begin{equation}
\tilde S=S^t, \quad S'=(S^{-1})^{t_1},\quad S''=(S^{t_2})^{-1},\label{eq3.12}
\end{equation}
where $t$ and $t_i$ mean the full and partial transpositions
respectively. Clearly, all the pentagon relations (\ref{eq3.2}) and
 (\ref{eq3.8}) are
reduced to (\ref{eq1.2}). Thus, formula (\ref{eq3.5})
enables to construct solutions to the Yang-Baxter relation from
 invertible and cross-invertible\footnote{ By cross-invertibility
  we mean the invertibility of the partially transposed matrix
 $S^{t_2}$} solutions to the pentagon relation (\ref{eq1.2}).

The above construction can be generalized also to the non-constant case.
For this consider some invertible and cross-invertible
 solution to the non-constant pentagon relation (see \cite{M} for the
 definition,
 and \cite{K}, for the example) which can be written in the following form:
\begin{eqnarray}
S_{12}(z_0,z_1,z_2,z_3)S_{13}(z_0,z_1,z_3,z_4)S_{23}(z_1,z_2,z_3,z_4)
&\nonumber\\
=S_{23}(z_0,z_2,z_3,z_4)S_{12}(z_0,z_1,z_2,z_4),&  \label{eq3.13}
\end{eqnarray}
where $z_0,\ldots,z_4$ are some parameters of any nature.
Substituting non-constant $S$'s into (\ref{eq3.5}), we obtain an $R$-matrix,
 depending of sixteen parameters:
\begin{eqnarray}
R_{12,34}(\hat z)=(S_{14}^{t_4}(z_{11},\ldots,z_{14}))^{-1}S_{13}(z_{21},
\ldots,z_{24})&\nonumber \\
\times S_{24}^t(z_{31},\ldots,z_{34})
(S_{23}^{-1}(z_{41},\ldots,z_{44}))^{t_2},& \label{eq3.14}
\end{eqnarray}
where $\hat z$ is a 4-by-4 matrix with entries $z_{ij}$ for $i,j=1,\ldots,4$.
Substituting six such $R$-matrices with different arguments into (\ref{eq3.7}),
we get the following equation:
\begin{equation}
R_{11',22'}(\hat x)R_{11',33'}(\hat y)R_{22',33'}(\hat z)
=R_{22',33'}(\hat z')R_{11',33'}(\hat y')R_{11',22'}(\hat x'), \label{eq3.15}
\end{equation}
which is satisfied provided the matrices $\hat x,\ldots, \hat x'$ are
constrained in the special way:
\begin{eqnarray}
&x_{23}=x_{43}=z_{12}=z_{22}=x'_{23}=x'_{43}=z'_{12}=z'_{22},\nonumber\\
&x_{13}=x_{33}=z_{32}=z_{42}=x'_{13}=x'_{33}=z'_{32}=z'_{42},\nonumber
\label{eq3.16}
\end{eqnarray}
\begin{eqnarray}
&&x_{22}=x_{42}=x'_{22}=x'_{42}=y_{i2},\quad
x_{12}=x_{32}=x'_{12}=x'_{32}=y'_{i2},\nonumber\\
&&z_{13}=z_{23}=z'_{13}=z'_{23}=y_{i3},\quad
z_{33}=z_{43}=z'_{33}=z'_{43}=y'_{i3},\quad i=1,\dots,4,\nonumber\label{eq3.17}
\end{eqnarray}
\begin{eqnarray}
&x_{11}=x'_{11}=y'_{11}=y'_{21}=z'_{31}=z'_{41},\quad
z_{14}=z'_{14}=x'_{24}=x'_{44}=y_{24}=y_{44},\nonumber\\
&x_{21}=x'_{21}=y_{11}=y_{21}=z'_{11}=z'_{21},\quad
z_{24}=z'_{24}=x_{24}=x_{44}=y_{14}=y_{34},\nonumber\\
&x_{31}=x'_{31}=y'_{31}=y'_{41}=z_{31}=z_{41},\quad
z_{34}=z'_{34}=x_{14}=x_{24}=y'_{14}=y'_{34},\nonumber\\
&x_{41}=x'_{41}=y_{31}=y_{41}=z_{11}=z_{21},\quad
z_{44}=z'_{44}=x'_{14}=x'_{24}=y'_{24}=y'_{44}.\label{eq3.18}
\end{eqnarray}
Apparently, there should exist a reparametrization, simplifying these
constraints, and admitting some ``rapidity'' picture.

To conclude the section note, that formula (\ref{eq3.5}) has a typical ``box''
structure, used for the construction of solutions to the Yang-Baxter
equation in terms of those to the ``twisted'' Yang-Baxter equations \cite{FM}.

\section{Examples}

\label{sec4}
\setcounter{equation}{0}
Consider several examples, realizing the general constructions of
Section~\ref{sec2}.

For a given finite group ${\cal G}$ consider it's group algebra as a
 Hopf algebra ${\cal A}$. The multiplication relations of the
corresponding Heisenberg double are as follows:
\begin{equation}
e_ge_h=e_{gh},\quad e^ge^h=\delta_{g,h}e^{g},\quad e^{g}e_h=e_he^{gh},\quad
g,h\in{\cal G}.\label{eq4.1}
\end{equation}
This example generates the infinite sequence of finite dimensional
solutions to the pentagon relation (\ref{eq1.2}). Particularly, in the adjoint
representation (\ref{eqadj}), the matrix elements of the canonical element
read:
\begin{equation}
\langle g,h|S|g',h'\rangle=\delta^{g'}_{gh}\delta^{h'}_h.
\end{equation}

Let Hopf algebra ${\cal A}$ be now the algebra ${\bf C}[x]$ with
the co-multiplication
\begin{equation}
\Delta(x)=x\otimes 1+1\otimes x.         \label{eq4.3}
\end{equation}
For a linear base take the normalized monomials:
\begin{equation}
e_m=x^m/m!,\quad m=0,1,\ldots    \label{eq4.4}
\end{equation}
The dual algebra ${\cal A^*}$ is also  the algebra  ${\bf C}[\overline x]$ with
 the co-multiplication of the form (\ref{eq4.3}), the dual base being
\begin{equation}
e^m=\overline x^m,\quad m=0,1,\ldots    \label{eq4.5}
\end{equation}
The Heisenberg double $H({\cal A})$ is defined by the Heisenberg
permutation relation between $x$ and $\overline x$:
\begin{equation}
x\overline x-\overline xx=1,                             \label{eq4.6}
\end{equation}
and the canonical element is just the usual exponent:
\begin{equation}
S=\sum_{m=0}^\infty e_m\otimes e^m=\exp(x\otimes\overline x).\label{eq4.7}
\end{equation}
The pentagon relation (\ref{eq1.2}) is an evident consequence of the
permutation relation (\ref{eq4.6}).

The next example is less trivial. Let algebra ${\cal A}$ be a deformed
universal enveloping algebra of the Lie algebra:
\begin{equation}
HE-EH=E,          \label{eq4.8}
\end{equation}
with the co-multiplications \cite{D}:
\begin{equation}
\Delta(H)=H \otimes 1+1\otimes H,\quad \Delta(E)=E\otimes\exp(hH)+1\otimes E,
                                   \label{eq4.9}
\end{equation}
where the complex parameter $h$ is chosen to have a positive real part.
 This algebra coincides with Borel subalgebra of $U_q(sl(2))$ where
$q=\exp(-h)$.
For a linear base take again normalized monomials:
\begin{equation}
e_{m,n}=H^mE^n/m!(q)_n, \quad m,n=0,1,\ldots,          \label{eq4.10}
\end{equation}
where
\begin{equation}
(q)_n=\cases{1&$n=0$;\cr
              (1-q)\ldots (1-q^n)&$n>0$.\cr}      \label{eq4.11}
\end{equation}
The dual algebra is generated also by two elements $\overline H$ and $F$ with
the permutation relation:
\begin{equation}
\overline HF-F\overline H=-hF,          \label{eq4.12}
\end{equation}
and co-multiplications:
\begin{equation}
\Delta(\overline H)=\overline H \otimes 1+1\otimes\overline H,\quad
\Delta(F)=F\otimes\exp(-\overline H)+1\otimes F,       \label{eq4.13}
\end{equation}
the dual base being
\begin{equation}
e^{m,n}=\overline H^mF^n,\quad m,n=0,1,\ldots  \label{eq4.14}
\end{equation}
The permutation relations of the Heisenberg double $H({\cal A})$ read
\begin{equation}
H\overline H-\overline HH=1,\quad E\overline H=\overline HE,\label{eq4.15}
\end{equation}
\begin{equation}
HF-FH=-F,\quad EF-FE=(1-q)q^{-H}.              \label{eq4.16}
\end{equation}
The canonical element is given by the formula:
\begin{equation}
S=\sum_{m,n=0}^\infty e_{m,n}\otimes e^{m,n}=\exp(H\otimes\overline H)
(E\otimes F;q)_\infty^{-1}, \label{eq4.17}
\end{equation}
where
\begin{equation}
(x;q)_\infty=(1-x)(1-xq)\ldots          \label{eq4.18}
\end{equation}
The pentagon relation (\ref{eq1.2}) for the element (\ref{eq4.17}) can be
 rewritten by the
use of above permutation relations in the following form:
\begin{equation}
(U;q)_\infty([U,V]/(1-q);q)_\infty (V;q)_\infty=(V;q)_\infty(U;q)_\infty,
                                                     \label{eq4.19}
\end{equation}
where the square brackets denote the commutator, and
\begin{equation}
U=E_2F_3,\quad V=E_1F_2.    \label{eq4.20}
\end{equation}
Operators $U$ and $V$ satisfy the following algebraic relations:
\begin{equation}
W\equiv UV-qVU,\quad [U,W]=[V,W]=0,      \label{eq4.21}
\end{equation}
which mean that the element $W$ lies in the center of the algebra,
generated by operators $U$ and $V$. In particular case, where $W=0$,
 formula (\ref{eq4.19}) coincides with the quantum dilogarithm identity of
 \cite{FK}. The generalized form (\ref{eq4.19}) of the latter
has been found earlier by Volkov \cite{V}.

\end{document}